\documentclass[12pt]{article}
\usepackage{amssymb}
\usepackage[left,pagewise,mathlines]{lineno}

\newcommand \dis{\displaystyle}

\title
{\bf The Forwarding Indices \\ of Graphs -- a Survey
\thanks {The work was
supported partially by NNSF of China (No. 11071233).} }
\author
{{Jun-Ming Xu\footnote{Corresponding author: xujm@ustc.edu.cn (J.-M.
Xu)}}\\
{\small School of Mathematical Sciences}  \\
{\small University of Science and Technology of China}   \\
{\small Hefei, Anhui, 230026, China}  \\
{\small E-mail addresses: xujm@ustc.edu.cn}\\ \\
{Min Xu
}\\
{\small School of Mathematical Sciences}  \\
{\small Beijing Normal University}   \\
{\small Beijing, 100875, China}  \\
{\small E-mail addresses: xum@bnu.edu.cn} }

\date{}

\begin{document}


\maketitle

\begin{abstract}

{\hskip20pt A routing $R$ of a given connected graph $G$ of order
$n$ is a collection of $n(n-1)$ simple paths connecting every
ordered pair of vertices of $G$. The vertex-forwarding index
$\xi(G,R)$ of $G$ with respect to $R$ is defined as the maximum
number of paths in $R$ passing through any vertex of $G$. The
vertex-forwarding index $\xi(G)$ of $G$ is defined as the minimum
$\xi(G,R)$ over all routing $R$'s of $G$. Similarly, the
edge-forwarding index $ \pi(G,R)$ of $G$ with respect to $R$ is the
maximum number of paths in $R$ passing through any edge of $G$. The
edge-forwarding index $\pi(G)$ of $G$ is the minimum $\pi(G,R)$ over
all routing $R$'s of $G$. The vertex-forwarding index or the
edge-forwarding index corresponds to the maximum load of the graph.
Therefore, it is important to find routings minimizing these indices
and thus has received much research attention in the past ten years
and more. In this paper we survey some known results on these
forwarding indices, further research problems and several
conjectures.}



\end{abstract}

\section{Introduction }
In a communication network, the message delivery system must find a
route along which to send each message from its source to its
destination. The time required to send a message along the fixed
route is approximately dominated by the message processing time at
either end-vertex, intermediate vertices on the fixed route relay
messages without doing any extensive processing. Metaphorically
speaking, the intermediate vertices pass on the message without
having to open its envelope. Thus, to a first approximation, the
time required to send a message along a fixed route is independent
of the length of the route. Such a simple forwarding function can be
built into fast special-purpose hardware, yielding the desired high
overall network performance.

For a fully connected network, this issue is trivial since every
pair of processors has direct communication in such a network.
However, in general, it is not this situation. The network designer
must specify a set of routes for each pair $(x,y)$ of vertices in
advance, indicating a fixed route which carries the data transmitted
from the message source $x$ to the destination $y$. Such a choice of
routes is called a routing.

We follow~\cite{x03} for graph-theoretical terminology and notation
not defined here. A graph $G=(V,E)$ always means a simple and
connected graph, where $V=V(G)$ is the vertex-set and $E=E(G)$ is
the edge-set of $G$. It is well known that the underlying topology
of a communication network can be modelled by a connected graph
$G=(V,E)$, where $V$ is the set of processors and $E$ is the set of
communication links in the network.

Let $G$ be a connected graph of order $n$. A routing $R$ in $G$ is a
set of $n(n-1)$ fixed paths for all ordered pairs $(x,y)$ of
vertices of $G$. The path $R(x,y)$ specified by $R$ carries the data
transmitted from the source $x$ to the destination $y$. A routing
$R$ in $G$ is said to be minimal, denoted by $R_m$, if each of the
paths specified by $R$ is shortest; $R$ is said to be symmetric or
bidirectional, if for all vertices $x$ and $y$, path $R(y,x)$ is the
reverse of the path $R(x,y)$ specified by $R$; $R$ is said to be
consistent, if for any two vertices $x$ and $y$, and for each vertex
$z$ belonging to the path $R(x,y)$ specified by $R$, the path
$R(x,y)$ is the concatenation of the paths $R(x,z)$ and $R(z,y)$.

It is possible that the fixed paths specified by a given routing $R$
going through some vertex are too many, which means that the routing
$R$ loads the vertex too much. Load of any vertex is limited by
capacity of the vertex, for otherwise it would affect efficiency of
transmission, even result in malfunction of the network.

It seems quite natural that a ``good" routing should not load any
vertex too much, in the sense that not too many paths specified by
the routing should go through it. In order to measure the load of
a vertex, Chung, Coffiman, Reiman and Simon~\cite{ccrs87} proposed
the notion of the forwarding index.

Let $G$ be a graph with a give routing $R$ and $x$ be a vertex of
$G$. The load of $x$ with respect to $R$, denoted by $\xi_x(G,R)$,
is defined as the number of the paths specified by $R$ going
through $x$. The parameter
$$
\xi(G,R)=\max\{\xi_x(G,R):\ x\in V(G)\}
$$
is called the forwarding index of $(G,R)$, and the parameter
$$
\xi(G)=\min\{\xi(G,R):\ \forall\ R\},\quad
\xi_m(G)=\min\{\xi(G,R_m):\ \forall\ R_m\}
$$
is called the forwarding index of $G$.

Similar problems are studied for edges by Heydemann, Meyer and
Sotteau~\cite{hms89}. The load of an edge $e$ with respect to $R$,
denoted by $\pi_e(G,R)$, is defined as the number of the paths
specified by $R$ which go through it. The edge-forwarding index of
$(G,R)$, denoted by $\pi(G,R)$, is the maximum number of paths
specified by $R$ going through any edge of $G$, i.e.,
$$
\pi(G,R)=\max\{\pi_e(G,R):\ e\in E(G)\};
$$
and the edge-forwarding index of $G$ is defined as
$$
\pi(G)=\min\{\pi(G,R):\ \forall\ R\},\quad
\pi_m(G)=\min\{\pi(G,R_m):\ \forall\ R_m\}.
$$

Clearly, $\xi(G)\le \xi_m(G)$ and $\pi(G)\le \xi_m(G)$. The
equality however does not always holds.

The original research of the forwarding indices is motivated by the
problem of maximizing network capacity~\cite{ccrs87}.  Maximizing
network capacity clearly reduces to minimizing vertex-forwarding
index or edge-forwarding index of a routing. Thus, whether or not
the network capacity could be fully used will depend on the choice
of a routing. Beyond a doubt, a ``good" routing should have a small
vertex-forwarding index and edge-forwarding index. Thus it becomes
very significant, theoretically and practically, to compute the
vertex-forwarding index and the edge-forwarding index of a given
graph and has received much attention the recent ten years and more.

Generally, computing the forwarding index of a graph is very
difficult. In this paper, we survey some known results on these
forwarding indices, further research problems, several conjectures,
difficulty and relations to other topics in graph theory.

Since forwarding indices were first defined for a graph, that is,
an undirected graph~\cite{ccrs87}, most of the results in the
literature are given for graphs instead of digraphs, but they can
be easily extended to digraphs. Nevertheless, we give here most of
the results for graphs as they appear in the literature.


\section{Basic Problems and Results}                                               

\noindent{\bf 2.1.\quad NP-completeness}\quad

\vskip6pt Chung, Coffiman, Reiman and Simon~\cite{ccrs87} asked
whether the problem of computing the forwarding index of a graph is
an NP-complete problem. Following~\cite{gj79}, we state this problem
as follows.

\vskip6pt {\bf Problem 2.1}\quad  Forwarding Index Problem,

\vskip6pt {\bf Instance}:\ A graph $G$ and an integer $k$.

\vskip6pt {\bf Question}:\ $\xi(G)\le k$?

\vskip6pt Heydemann, Meyer, Sotteau and Opatrn\'y~\cite{hmso94}
first showed that Problem 2.1 is NP-complete for graphs of
diameter at least $4$ when the routings considered are restricted
shortest, consistent and symmetric; a P-problem for graphs of
diameter $2$ when the routings considered are restricted to be
shortest. Saad~\cite{s93} proved that Problem 2.1 is NP-complete for
for general routings even if the diameter of the graph is $2$.
However, Problem 2.1 has not yet been solved for graphs of
$3$ when the routings considered are restricted to be shortest
be minimal and/or, consistent and/or symmetric.

\vskip6pt The same problem was also suggested by Heydemann, Meyer
and Sotteau~\cite{hms89}.

\vskip6pt {\bf Problem 2.2}\quad Edge-Forwarding Index Problem,

\vskip6pt {\bf Instance}:\ A graph $G$ and an integer $k$.

\vskip6pt {\bf Question}:\ $\pi(G)\le k$?

\vskip6pt Heydemann, Meyer, Sotteau and Opatrn\'y~\cite{hmso94}
showed that Problem 2.2 is NP-complete for graphs of diameter at
least $3$ when the routings considered are restricted to be
minimal, consistent and symmetric; a P-problem for graphs of
diameter $2$ when the routings considered are restricted to be
minimal.

\vskip20pt \noindent{\bf 2.2.\quad Basic Bounds and Relations}\quad

\vskip10pt For a given connected graph $G$ of order $n$, set
$$
A(G)=\frac {1}{n}\sum\limits_{u\in V} \left(\sum\limits_{v\in
V\setminus\{u\}}(d_{G}(u,v)-1)\right),
$$
and
$$
B(G)=\frac{1}{\varepsilon}\sum\limits_{(u,\,v)\in V\times
V}d_{G}(u,v).
$$

The following bounds of $\xi (G)$ and $\pi(G)$ were first
established by Chung, Coffiman, Reiman and Simon~\cite{ccrs87} and
Heydemann, Meyer and Sotteau~\cite{hms89}, respectively.

\vskip6pt {\bf Theorem 2.3}\ (Chung {\it et al}~\cite{ccrs87})
\quad Let $G$ be a connected graph of order $n$. Then
\begin{equation}
\label{e1}                                                           
A(G) \le \xi(G)\le (n-1)(n-2),
\end{equation}
and the equality $\xi_{G}=\xi_m(G)=A(G)$ is true if and only if
there exists a minimal routing in $G$ which induces the same load on
every vertex. The graph that attains this upper bound is a star
$K_{1,n-1}$.

\vskip6pt {\bf Theorem 2.4}\ (Heydemann {\it et al}~\cite{hms89})
\quad Let $G$ be a connected graph of order $n$. Then
\begin{equation}
\label{e2}                                                           
B(G)\le \pi(G) \le \left\lfloor \frac 12n^2\right\rfloor,
\end{equation}
and the equality $\pi(G)=\pi_m(G)=B(G)$ is true if and only if there
exists a minimal routing in $G$ which induces the same load on every
edge. The graph that attains this upper bound is a complete
bipartite graph $K_{\frac n2, \frac n2}$.

Recently, Xu {\it et al}.~\cite{xzdy07} have showed the star
$K_{1,n-1}$ is a unique graph that attains the upper bound in
(\ref{e1}).

\vskip6pt {\bf Problem 2.5}\quad Note that the upper bound given in
(\ref{e2}) can be attained. Give a characterization of graphs whose
vertex- or edge-forwarding indices attain the upper bound in
(\ref{e2}).

\vskip6pt Although the two concepts of vertex- and edge-forwarding
index are similar, no interesting relationships is known between
them except the following trivial inequalities.

\vskip6pt {\bf Theorem 2.6}\quad (Heydemann {\it et
al}~\cite{hms89}) For any connected undirected graph $G$ of order
$n$, maximum degree $\Delta$, minimum degree $\delta$,

(a)\ $2\xi(G)+2(n-1)\le \Delta\pi(G)$;

(b)\ $\pi(G)\le \xi(G)+2(n-1)$;

(c)\ $\pi_m(G)\le \xi_m(G)+2(n-\delta)$.

\noindent The inequality in (a) is also valid for minimal
routings.

\vskip6pt No nontrivial graph is found for which the forwarding
indices hold one of the above equalities. Thus, it is necessary to
further investigate the relationships between $\pi(G)$ and
$\xi(G)$ or between $\pi_m(G)$ and $\xi_m(G)$.

\vskip6pt {\bf Problem 2.7}\quad For a graph $G$ and its line graph
$L(G)$, investigate the relationships between $\xi(G)$ and
$\pi(L(G))$ or between $\xi_m(G)$ and $\pi_m(L(G))$.

\vskip20pt \noindent{\bf 2.3.\quad Optimal Graphs}\quad

\vskip10pt A graph $G$ is said to be vertex-optimal if
$\pi(G)=A(G)$, and edge-optimal if $\xi(G)=B(G)$. Note that if $R_m$
is a routing of $G$ such that $\pi(G,R_m)=A(G)$, then
\begin{equation}
\label{e3}                                                         
\xi(G)=\sum\limits_{y\in V}d(G;x,y)-(n-1), \quad \forall \ x\in V.
\end{equation}
Heydemann {\it et al}~\cite{hms89} showed that the equality (3) is
valid for any Cayley graph.

\vskip6pt {\bf Theorem 2.8}\quad Let $G$ be a connected Cayley graph
with order $n$. Then
\begin{equation}
\label{e4}                                                           
\xi (G)=\xi_m(G)=\sum\limits_{y\in V}d(G;x,y)-(n-1), \quad
\forall\ x\in V.
\end{equation}

From Theorem 2.8, Cayley graphs are vertex-optimal. Some results and
problems on the forwarding indices of vertex-transitive or Cayley
graphs, an excellent survey on this subject has been given by
Heydemann~\cite{h97}.

Heydemann {\it et al}~\cite{hmso94} have constructed a class of
graphs for which the vertex-forwarding index is not given by a
minimal consistent routing. Thus, they suggested the following
problems worthy of being considered.

\vskip6pt {\bf Problem 2.9}\quad (Heydemann {\it et
al}~\cite{hmso94}) For which graph or digraph $G$ does there exist a
minimal consistent routing $R$ such that $\xi_m(G)=\xi(G,R)$ or a
consistent routing $R$ such that $\xi(G)=\xi(G,R)$?

\vskip6pt Heydemann et al~\cite{hmso94} have ever conjectured that
in any vertex-transitive graph $G$, there exists a minimal routing
$R_m$ in which the equality (4) holds.

The conjecture has attracted many researchers for ten years and
more without a complete success until 2002. Shim, \v{S}ir\'a\v{n}
and \v{Z}erovnik~\cite{ssz02} disproved this conjecture by
constructing an infinite family of counterexamples, that is,
$K_p\oplus P(10,2)$ for any $q\not\equiv0$\,(mod $3$), where
$P(10,2)$ is the generalized Petersen graph and the symbol
$\oplus$ denotes the strong product.

Gauyacq~\cite{g971, g97, gmr98} defined a class of quasi-Cayley
graphs, a new class of vertex-transitive graphs, which contain
Cayley graphs, and are vertex-optimal. Sol\'e~\cite{so94}
constructed an infinite family of graphs, the so-called orbital
regular graphs, which are edge-optimal. We state the results of
Gauyacq and Sol\'e as the following theorem.

\vskip6pt {\bf Theorem 2.10}\quad Any quasi-Cayley graph is
vertex-optimal, and any orbital regular graph is edge-optimal.

\vskip6pt However, we have not yet known whether a quasi-Cayley
graph is edge-optimal and not known whether an orbital regular
graph is vertex-optimal. Thus, we suggest to investigate the
following problem.

\vskip6pt {\bf Problem 2.11}\quad Investigate whether a quasi-Cayley
graph is edge-optimal and an orbital regular graph is
vertex-optimal.

\vskip6pt Considering $\pi(K_2\times K_p)$ for $p\ge 3$, Heydemann
{\it et al}~\cite{hms89} found that the equality (4) is not valid
for $\pi(G)$, and proposed the following conjecture.

\vskip6pt {\bf Conjecture 2.12}\ (Heydemann et al~\cite{hms89})\ For
any distance-transitive graph $G=(V,E)$, there exists a minimal
routing $R_m$ for which,
$$ \pi(G)=\pi (G,R_m)=\left\lceil\frac
{n}{\varepsilon}\sum\limits_{y\in V}d(G;x,y)\right\rceil, \quad
\forall\ x\in V.
$$

{\bf Conjecture 2.13}\quad (Heydemann {\it et al}~\cite{hms89}) For
any distance-transitive graph $G=(V,E)$, there exists a minimal
routing in which we have both

(a)\ the load of all vertices is the same, and then,
$$
\xi(G)=\xi_m(G)=\sum\limits_{y\in V}d(G;x,y)-(n-1), \quad \forall\
x\in V.
$$

(b)\ the load of all the edges is almost the same (difference of
at most one) and then,
$$
\pi(G)=\pi_m(G)=\left\lceil\frac {n}{\varepsilon}\sum\limits_{y\in
V}d(G;x,y)\right\rceil, \quad \forall\ x\in V.
$$

\vskip20pt \noindent{\bf 2.4.\quad For Cartesian Product
Graphs}\quad

\vskip10pt The cartesian product can preserve many desirable
properties of the factor graphs. A number of important
graph-theoretic parameters, such as degree, diameter and
connectivity, can be easily calculated from the factor graphs. In
particular, the cartesian product of vertex-transitive (resp.
Cayley) graphs is still a vertex-transitive (resp. Cayley) graph
(see Section 2.3 in \cite{x01}). Since quasi-Cayley graphs are
vertex-transitive, the cartesian product of quasi-Cayley graphs is
still a quasi-Cayley graph. Thus, determining the forwarding
indices of the cartesian product graphs is of interest. Heydemann
{\it et al}~\cite{hms89} obtained the following results first.

\vskip6pt {\bf Theorem 2.14}\quad Let $G$ and $G'$ be two
connected graphs with order $n$ and $n'$, respectively. Then

(a)\ $\xi(G\times G')\le n\ \xi(G')+n'\ \xi(G)+(n-1)(n'-1)$;

(b)\ $\pi(G\times G')\le \max\{n\pi(G'),\ n'\pi(G)\}.$

The inequalities are also valid for minimal routings. Moreover,
the equality in (a) holds if both $G$ and $G'$ are Cayley
digraphs.

\vskip6pt Recently, Xu {\it et al}~\cite{xxh06} have considered the
cartesian product of $k$ graphs and obtained the following results.

\vskip6pt {\bf Theorem 2.15}\quad Let $G=G_1\times G_2\times
\cdots \times G_k$. Then

(a)\ $G$ is vertex-optimal and if $G_i$ is vertex-optimal for every
$i=1,2,\cdots,k$ then
 $$\xi(G)= \sum\limits_{i=1}^k n_1n_2\cdots
n_{i-1}(\xi_i-1)n_{i+1}\cdots n_k+(k-1)n_1n_2\cdots n_k+1;
 $$

(b)\ $G$ is edge-optimal and if $G_i$ is edge-optimal for every
$i=1,2,\cdots,k$ then
 $$
 \pi(G)= \max\limits_{1\le i\le k}
\{n_1n_2\cdots n_{i-1}\pi_in_{i+1}\cdots n_k\}.
$$

\vskip6pt By Theorem 2.15 and Theorem 2.10, the cartesian product of
quasi-Cayley graphs is vertex-optimal and the cartesian product of
orbital regular graphs is edge-optimal.

\section{Connectivity Constraint}

In this section, we survey the known results of the forwarding
indices of $k$-connected or $k$-edge-connected graphs.

\vskip20pt \noindent{\bf 3.1.\quad $\kappa$-connected Graphs}\quad

\vskip10pt{\bf Theorem 3.1}\quad If $G$ is a $2$-connected graph
of order $n$, then

(a)\ $\xi (G)\le \frac 12\,(n-2)(n-3)$, this bound is best
possible in view of $K_{2,n-2}$ (Heydemann {\it et
al}~\cite{hms89});

(b)\ $\xi_m(G)\le n^2-7n+12$ for $n\ge 6$ and diameter $2$, this
bound is best possible since it is reached for a wheel of order
$n$ minus one edge with both ends of degree $3$ (Heydemann {\it et
al}~\cite{hms89});

(c)\ $\xi_m(G)\le n^2-7n+12$ for $n\ge 7$ (Heydemann {\it et
al}~\cite{hmos92});

(d)\ $\pi(G)\le \lfloor\frac 14\,n^2\rfloor$ and this bound is
best possible in view of the cycle $C_n$ (Heydemann {\it et
al}~\cite{hmos92}).

\vskip6pt{\bf  Theorem 3.2}\quad (Heydemann {\it et
al}~\cite{hmos92}) $\xi(G)\le n^2-(2k+1)n+2k$ for any
$k$-connected graph $G$ of order $n$ with $k\ge 3$ and $n\ge
8k-10$.

\vskip6pt Heydemann {\it et al}~\cite{hmos92} proposed the
following research problem.

\vskip6pt {\bf Problem 3.3}\quad Find the best upper bound
$f(n,k), g(n,k), h(n,k)$ and $s(n,k)$ such that for any
$k$-connected graph $G$ of order $n$ with $k\ge 2$,\ $\xi(G)\le
f(n,k), \xi_m(G)$ $\le g(n,k), \pi(G)\le h(n,k)$ and $\pi_m(G)\le
s(n,k)$ for $n$ large enough compared to $k$.

\vskip6pt {\bf Theorem 3.4}\quad (de la Vega and
Manoussakis~\cite{fm92}) For any integer $k\ge 1$,

(a)\ $f(n,k)\le (n-1) \left\lceil\frac 1k\,(n-k-1)\right\rceil$;

(b)\ $g(n,k)\le \frac 12n^2-(k-1)n+\frac 38(k-1)^2$ if $n$ is
substantially larger than $k$;

(c)\ $h(n,k)\le n\left\lceil\frac 1k\, (n-k-1)\right\rceil$.

\vskip6pt Recently, Zhou {\it et al}.~\cite{zxx08} have improved the
upper bounds of $f(n,k)$ and $h(n,k)$ in Theorem 3.4 as follows.

\vskip6pt{\bf Theorem 3.5} \ If $G$ is a $k$-connected graph of
order $n$ with the maximum degree $\Delta$, then $\xi (G)\le
(n-1)\lceil(n-k-1)/k\rceil-(n-\Delta-1)$ and $\pi (G)\le
n\lceil(n-k-1)/k\rceil-(n-\Delta)$.

\vskip6pt {\bf Conjecture 3.6}\quad (de la Vega and
Manoussakis~\cite{fm92}) For any positive integer $k$,

(a)\ $f(n,k)\le\lceil\frac 1k\,(n-k)(n-k-1)\rceil$ for $n\ge 2k\ge
2$, which would be best possible in view of the complete bipartite
graph $K_{k,n-k}$;

(b)\ there exists a function $q(k)$ such that if $n\ge q(k)$, then
$g(n,k)\le \frac 12n^2-(k-1)n-\frac 32k^2+k+\frac 72$;

(c)\ $h(n,k)\le \lceil\frac {n^2}{2k}\rceil$ for $n\ge 2k\ge 2$,
which would be best possible in view of the graph obtained from
two complete graphs $K_m$ plus a matching $e_1,e_2,\cdots,e_k$
between them, $m\ge k$.

\vskip6pt It can be easily verified that the conjecture (a) and (c)
are true for $k=1$ and $k=2$. Recently, Zhou {\it et
al}.~\cite{zxx08} have proved that Conjecture 3.6 (a) is true for
$k=3$, that is,

\vskip6pt{\bf Theorem 3.7}\ If $G$ is a $3$-regular and
$3$-connected graph of order $n\ge 4$. Then $\xi(G)\le \lceil(n-3)(n-4)/3\rceil$. 

\vskip20pt \noindent{\bf 3.2.\quad $\lambda$-edge-connected
Graphs}\quad

\vskip10pt{\bf Theorem 3.8}\quad If $G$ is a $2$-edge-connected
graph of order $n$, then

(a)\ $\pi_m(G)\le \lfloor\frac 12\,n^2-n+\frac 12\rfloor$
(Heydemann {\it et al}~\cite{hms89,  hmos92});

(b)\ $\pi(G)\le \lfloor\frac 14\,n^2\rfloor$ (Cai~\cite{cai90}).

\vskip6pt Heydemann {\it et al}~\cite{hms89} conjectured that for
any $\lambda$-edge-connected graph $G$ of order $n$, $\pi(G)\le
\lfloor\frac 12\,n^2-(\lambda-1)n+\frac 12\, (\lambda-1)^2
\rfloor$. Latter, Heydemann, Meyer, Opatrn\'y and
Sotteau~\cite{hmos92} gave a counterexample and proposed the
following conjecture.

\vskip6pt{\bf  Conjecture 3.9}\quad for any $\lambda$-edge-connected
graph $G$ of order $n$ with $\lambda\ge 3$ and $n\ge 3\lambda+3$,
$$
\pi_m(G)=\max\left\{\left\lceil\frac
{n^2}2\right\rceil-n-2(\lambda-1)^2,\ \left\lfloor\frac
{n^2}2\right\rfloor-2n+5\lambda-\frac 32(\lambda^2+1)\right\}.
$$

The same problem as the ones in Problem 3.3 can be considered for
$\lambda$-edge-connected graphs.

\vskip6pt {\bf Problem 3.10}\quad Find the best upper bound
$f'(n,\lambda), g'(n,\lambda), h'(n,\lambda)$ and $s'(n,\lambda)$
such that for any $\lambda$-connected graph $G$ of order $n$ with
$\lambda\ge 2$,\ $\xi(G)\le f'(n,\lambda),\xi_m(G)$ $\le
g'(n,\lambda), \pi(G)\le h'(n,\lambda)$ and $\pi_m(G)\le
s'(n,\lambda)$ for $n$ large enough compared to $\lambda$.

\vskip6pt The following theorem is the only result we have known
as far on this problem.

\vskip6pt {\bf Theorem 3.11}\quad (de la Vega and
Manoussakis~\cite{fm92}) For any integer $\lambda\ge 3$,
$$
g'(n,\lambda)=\left\lceil\frac
{n^2}2\right\rceil-n-2(\lambda-1)^2\ {\rm for}\
n\ge\max\left\{3\lambda+3,\ \frac 12(\lambda+1)^2\right\}.
$$

\vskip20pt \noindent{\bf 3.3.\quad Strongly Connected Digraphs}\quad

\vskip10pt It is clear that the notion of the forwarding indices
can be similarly defined for digraphs. Many general results, such
as Theorem 2.3 and Theorem 2.4 are valid for digraphs. Manoussakis
and Tuza~\cite{mt962} consider the forwarding index of strongly
$k$-connected digraphs and obtained the following result similar
to Theorem 2.4.

\vskip6pt {\bf Theorem 3.12}\quad Let $D=(V,E)$ be a strong digraph
of order $n$. Then

(1) $B(D)\le \pi(D)\le \pi_m(D)\le (n-1)(n-2)+1$, and

(2) The equalities $\pi(D)=\pi_m(D)=B(D)$ are true if and only if
there exists a minimal routing in $D$ which induces the same load
on every edge.

\vskip6pt In addition to validity of Theorem 3.2 for digraphs,
they obtained the following results.

\vskip6pt {\bf Theorem 3.13}\quad Let $D$ be a $k$-connected digraph
of order $n\ge 3$, and $k\ge 1$. Then

(a)\ $\pi (D)\le (n-1)\lceil\frac 1k(n-k-1)\rceil+1$;

(b)\ $\xi_m(D)\le n^2-(2k+1)n+2k$ for $n\ge 2k+1$;

(c)\ $\pi_m(D)\le n^2-(3K+2)n+4k+3$ for $n\ge 4k-1$.

\section{Degree Constraint}

Although Saad~\cite{s93} showed that for any graph determining the
forwarding index problem is NP-complete, yet many authors are
interested in the forwarding indices of a graph. Specially, it is
still of interest to determine the exact value of the forwarding
index with some graph-theoretical parameters. For example, Chung
{\it et al}~\cite{ccrs87}, Bouabdallah and Sotteau~\cite{bs93}
proposed to determine the minimum forwarding indices of
$(n,\Delta)$-graphs that has order $n$ and maximum degree at most
$\Delta$. Given $\Delta$ and $n$, let
$$
\begin{array}{rl}
 &\xi_{\Delta,n}=\min\{\xi(G):\ |V(G)|=n, \Delta(G)=\Delta\},\\
  &\pi_{\Delta,n}=\min\{\pi(G):\ |V(G)|=n, \Delta(G)=\Delta\}.
  \end{array}
$$

\vskip6pt \noindent{\bf 4.1.\quad Problems and Trivial Cases}

\vskip10pt {\bf Problem 4.1}\quad (Chung {\it et
al}~\cite{ccrs87}) Given $\Delta\ge 2$ and $n\ge 4$, determine
$\xi_{\Delta,n}$, and exhibit an $(n,\Delta)$-graph $G$ and $R$ of
$G$ for which $\xi(G,R)=\xi_{\Delta,n}$.

\vskip6pt {\bf Problem 4.2}\quad (Bouabdallah and
Sotteau~\cite{bs93}) Given $\Delta\ge 2$ and $n\ge 4$, determine
$\pi_{\Delta,n}$, and exhibit an $(n,\Delta)$-graph $G$ and $R$ of
$G$ for which $\pi(G,R)=\pi_{\Delta,n}$.

\vskip6pt For $\Delta \ge n-1$, we can fully connect a graph,
i.e., $G$ is a complete graph. In this case any routing $R$ can be
composed only of single-edge paths so that the minimum $\xi=0$ and
$\pi=2$ is achieved, that is, $\xi_{\Delta,n}=\xi(G,R)=0$ and
$\pi(G,R)=\pi_{\Delta,n}=2$ for $\Delta\ge n-1$.

For $\Delta =2$ the only connected graph fully utilizing the
degree constraint is easily seen to be a cycle. Because of the
simplicity of cycles, the for vertex-forwarding index problem can
be solve completely for $\Delta =2$.

\vskip6pt {\bf Theorem 4.3}\quad For all $n\ge 3$,\

(a)\ $\xi_{2,n}=\xi(C_n,R_m)=\left\lfloor\frac
14(n-1)^2\right\rfloor$;

(b)\ $\pi_{2,n}=\pi(C_n,R_m)=\left\lfloor\frac
14n^2\right\rfloor$.


\vskip20pt \noindent{\bf 4.2.\quad Results on $\xi_{\Delta,\,n}$}
\ $\Delta\ge 3$ \quad

\vskip10pt Problem 4.1 was solved for $n\le 15$ or any $n$ and
$\Delta$ with $\frac 13\,(n+4)\le \Delta\le n-1$ by Heydemann et
al~\cite{hms88}.

\vskip6pt {\bf Theorem 4.4}\quad (Heydemann {\it et
al}~\cite{hms88})

(a)\ if $n$ is even or $n$ odd and $\Delta $ even,
$\xi_{\Delta,\,n}=n-1-\Delta$ for $\Delta\ge \frac 13\,(n+1)$ or
for $n=12$ or $13$ and $\Delta =4$;

(b)\ if $n$ and $\Delta $ are odd, $\xi_{\Delta,\,n}=n-\Delta$ for
$\Delta\ge \frac 13\,(n+4)$ or for $n=13$ and $\Delta =5$.

\vskip6pt Problem 4.1 has not been completely solved for
$\Delta<\frac 13\,(n+4)$.

\vskip6pt {\bf Theorem 4.5}\quad (Heydemann {\it et
al}~\cite{hms88}) For any $n$ and $\Delta$,


(a)\ $\xi_{n-2p-1,\,n}=2p$ for any $n$ and $p$ such that $n\ge
3p+2$;

(b)\ $\xi_{2p+1,\,n}=n-2p-1$ for any odd $n$ such that $2p+1\le
n\le 6p-1$;

(c)\ $\xi_{2p,\,n}=n-2p-1$ for any $n$ and $p$ such that $2p+1\le
n\le 6p-1$ and $p\ge 3$;

(d)\ $\xi_{\Delta,\,n}\ge n-1-\Delta$;

(e)\ $\xi_{\Delta,\,n}=n-1-\Delta\Longrightarrow$ every
$(n,\Delta)$-graph $G$ such that $\xi(G)=\xi_{\Delta,\,n}$ is
$\Delta$-regular and diameter $2$;

(f)\ $\xi_{\Delta,\,n}\ge n-\Delta$ if $n$ and $\Delta$ are odd.

\vskip6pt An asymptotic result on $\xi_{\Delta,\,n}$ has been
given by Chung {\it et al}~\cite{ccrs87}.

\vskip6pt {\bf Theorem 4.6}\quad For any given $\Delta \ge 3$,
$$
[1+o(1)]n\log_{\Delta-1}n\le \xi_{\Delta,\,n}\le
\left[3+O\left(\dis\frac
1{\log\Delta}\right)\right]n\log_{\Delta}n,
$$
where the upper bound holds for $\Delta\ge 6$.

\vskip20pt \noindent{\bf 4.3.\quad Results on $\pi_{\Delta,\,n}$}
\ $\Delta\ge 3$ \quad

\vskip10pt Similar to Theorem 4.5, Bouabdallah and
Sotteau~\cite{bs93} obtained the following result on
$\pi_{\Delta,\,n}$

\vskip6pt {\bf Theorem 4.7}\quad For any $n$ and $\Delta\ge 3$,

(a)\ $\pi_{\Delta,\,n}\ge \lceil\frac {4(n-1)}{\Delta}\rceil-2$;

(b)\ $\pi_{\Delta,\,n}\ge \lceil\frac
{4(n-1)}{\Delta}\rceil-2\Longrightarrow$ every $(n,\Delta)$-graph
$G$ such that $\pi(G)=\pi_{\Delta,\,n}$ is $\Delta$-regular and
diameter $2$ and $G$ has a minimal routing for which the load of
all edges is the same;

(c)\ $\pi_{\Delta,\,n}\ge \lceil\frac {4n-2)}{\Delta}\rceil-2$ if
$n$ and $\Delta$ are odd;

(d)\ $\pi_{\Delta,\,n}\le \pi_{\Delta',\,n}$ for any $n$ and
$\Delta'\le n-1$ with $\Delta'<\Delta$.

\vskip6pt Problem 4.3 was solved for $n\le 15$ by Bouabdallah and
Sotteau~\cite{bs93}, who also obtained $\pi_{n-2, n}=3$ for any
$n\ge 6$, $n\ne 7$ and $\pi_{n-2, n}=4$ for any $n=4,5,7$. Recently,
Xu {\it et al}~\cite{xhx04} have determined $\pi_{n-2p-1,\,n}=8$ if
$3p+\lceil\frac 13\, p \rceil +1 \leq n \leq 3p+\lceil\frac 35\,p
\rceil$ and $\ge 2$. The authors in the two \cite{bs93} and
\cite{xhx04} have completely determined $\pi_{n-2p-1, n}$ for $n\ge
4p$ and $p\ge1$ except a little gap.

\vskip6pt {\bf Theorem 4.8}\quad  For any $p\geq 1$,
$$
\pi_{n-2p-1, n}=\left\{
\begin{array}{ll}
3, \ & {\rm if}\ n\geq 10p+1;\\
4, \ & {\rm if}\ 6p+1 \leq n < 10p+1;\\
5, \ & {\rm if}\ 4p+2\lceil\frac 13\, p \rceil +1 \leq n \leq 6p\,;\\
6, \ & {\rm if}\ 4p+1 \leq n \leq 4p+\lceil \frac 23\,p\rceil.
\end{array}\right.
$$

Note the value of $\pi_{n-2p-1,\,n}$ has not been determined for
$4p+\lceil \frac 23\,p \rceil+1\le n \le 4p+2\lceil\frac 13\,p
\rceil$. However, these two numbers are different only when
$p=3k+1$. Thus, we proposed the following conjecture.

\vskip6pt{\bf Conjecture 4.9}\quad For any $p\geq 1$,
$\pi_{n-2p-1, n}=5$ if $4p+\lceil\frac 23\,p \rceil+1\le n \le
4p+2\lceil\frac 13\, p \rceil$.

\vskip6pt {\bf Theorem 4.10}\ (Xu {\it et al}~\cite{xxh05})\quad For
any $p\geq 1$, we have
$$
\pi_{n-2p,\,n}=\left\{
\begin{array}{ll}
3, \ & {\rm if}\ n\geq 10p-2\ {\rm or}\ n=10p-4;\\
4, \ & {\rm if}\ 6p+1 \leq n < 10p-4\ or\ n=10p-3;\\
6, \ & {\rm if}\ 4p+1 \leq n \leq 4p+\lceil \frac 13\,(2p-1)
\rceil-2.
\end{array}\right.
$$

\vskip6pt An asymptotic result on $\pi_{\Delta,\,n}$ has been
given by Heydemann {\it et al}~\cite{hms89}.

\vskip6pt {\bf Theorem 4.11}\quad For any given $\Delta \ge 3$,
$$
\left[\frac 2{\Delta}+o(1)\right]\,\,n\,\log_{\Delta-1}n\le
\pi_{\Delta,\,n}\le 24\frac{\log_2(\Delta-1)}{\Delta}\, n
\log_{\Delta-1}n,
$$
where the upper bound holds for $\Delta\ge 6$.

\vskip20pt \noindent{\bf 4.4.\quad General Results Subject to Degree
and Diameter}\quad

\vskip6pt{\bf Theorem 4.12}\quad (Xu {\it et al}.~\cite{xzdy07}) For
any connected graph $G$ of order $n$ and maximum degree $\Delta$,
 $$
 \xi(G)\le (n-1)(n-2)-\left(2n-2-\Delta\left\lfloor 1+
\frac{n-1}{\Delta}\right\rfloor\right)\left\lfloor\frac{n-1}{\Delta}
\right\rfloor.
 $$

Considering a special case of $\Delta=n-1$ in Theorem 4.12, we
obtain the upper bound in (\ref{e1}) immediately.

\vskip6pt{\bf Theorem 4.13}\quad (Heydemann {\it et
al}~\cite{hms89}) Let $G$ be a graph of order $n$, maximum degree
$\Delta$ and diameter $d$,

(a)\ $\xi(G)\le \xi_m(G)\le (n-1)(n-2)-2(\varepsilon(G)-\Delta)$,

(b)\ $\xi(G)\le \xi_m(G)\le n^2-3n-\lfloor\frac
12\,d\rfloor^2-\lceil\frac 12\,d\rceil^2+d+2$.

\vskip6pt{\bf Theorem 4.14}\quad (Heydemann {\it et
al}~\cite{hms89}) If $G$ is a graph of order $n$ and diameter $2$
with no end vertex, then $\pi_m(G)\le 2n-4$.

\vskip6pt Manoussakis and Tuza~\cite{mt962} obtained some upper
bounds on the forwarding indies for digraphs subject to degree
constraints.

\vskip6pt{\bf Theorem 4.15}\quad Let $D$ be a strongly connected
digraph of order $n$ and minimum degree $\delta$. Then

(a)\ $\xi_m(D)\le n^2-(\delta+2)n+\delta+1$;

(b)\ $\pi_m(D)\le \max\{n^2-3n\delta+2\delta^2+\delta,\
n^2-(2\delta+3)n+\delta^2+4\delta+3\}$ if $n$ is sufficient large
compared to $\delta$.

Considering the minimum degree $\delta$ rather the maximum degree
$\Delta$, we can propose an analogy of $\xi_{\Delta,n}$ and
$\pi_{\Delta,n}$ as follows. Given $\Delta$ and $n$, let
$$
\begin{array}{rl}
 &\xi_{\delta,n}=\min\{\xi(G):\ |V(G)|=n, \delta(G)=\delta\},\\
  &\pi_{\delta,n}=\min\{\pi(G):\ |V(G)|=n, \delta(G)=\delta\}.
  \end{array}
$$

However, the problem determining $\xi_{\delta,n}$ and
$\pi_{\delta,n}$ is simple.

\vskip6pt {\bf Theorem 4.16} \quad (Xu {\it et al}.~\cite{xzdy07})\
For any $n$ and $\delta$ with $n>\delta\ge 1$,
 $$
 \xi_{\delta,n}=\left\lceil\frac{2(n-1-\delta)}{\delta}\right\rceil\quad
 {\rm and}\
\pi_{\delta,n}=\left\lceil\frac{2(n-1)}{\delta}\right\rceil.
 $$




\vskip10pt

\section{Difficulty, Methods and  Relations to Other Topics}

\vskip10pt \noindent{\bf 5.1.\quad Difficulty of Determining the
Forwarding Indices}\quad

As we have stated in Subsection 2.1, the problem of computing the
forwarding indices of a general graph is an NP-complete problem.
Also, for a given graph $G$, determining its forwarding indices
$\xi(G)$ and $\pi(G)$ is also very difficult.

The first difficulty is designing a routing $R$ such that
$\xi_x(G,R)$ for any $x\in V(G)$ or $\pi_x(G,R)$ for any $e\in E(G)$
can be conveniently computed. An ideal routing can be found by the
current algorithms for finding shortest paths. However, in general,
it is not always the case that the forwarding indices of a graph can
be obtained by a minimum routing.

For example, consider the wheel $W_7$ of order seven. The hub $x$,
other vertices $0,1,\cdots, 5$. A minimum and bidirectional routing
$R_m$ is defined as follows.
$$\left\{
\begin{array}{ll}
 R_m(i,i+2)=R(i+2,i)=(i,i+1,i+2)({\rm mod}\ 6), & i=0,1,\cdots,5;\\
 R_m(i,i+3)=R(i+3,i)=(i,x,i+3)({\rm mod}\ 6), & i=0,1,2;\\
 {\rm direct\ edge}, & {\rm otherwise}.
 \end{array}\right.
$$
Then,
$$
\xi_x(W_7,R_m)=6,\ \xi_i(W_7,R_m)=2,\ i=0,1,\cdots,5.
$$
Thus, we have $\xi(W_7,R_m)=6$.

However, if we define a routing $R$ that is the same as the minimum
routing $R_m$ except for $R(2,5)=(2,1,0,5),\ R(5,2)=(5,4,3,2)$. Then
the routing $R$ is not minimum. We have
$$
\xi_x(W_7,R)=4,\ \tau_2(W_7,R)=\tau_5(W_7,R)=2,\quad
\xi_i(W_7,R)=3,\ i=0,1,3,4.
$$
Thus, we have $\xi(W_7,R)=4<6=\xi(W_7,R_m)$.

The second difficulty is that the forwarding indices are always
attained by a bidirectional routing. For example, for the hypercube
$Q_n$ $(n\ge 2)$, $\xi(Q_n)=n2^{n-1}-(2^n-1)=2^{n-1}(n-2)+1$. Since
$2^{n-1}(n-2)+1$ is odd, $\xi(Q_n)$ can not be attained by a
bidirectional routing.

\vskip10pt \noindent{\bf 5.2.\quad Methods of Determining the
Forwarding Indices}\quad

To the knowledge of the author, one of the actual methods of
determining forwarding index is to compute the sum of all pairs of
vertices according to (\ref{e4}) for some Cayley graphs. In fact,
the forwarding indices of many Cayley graphs are determined by using
(\ref{e4}), for example, the folded cube~\cite{hxx05}, the augmented
cube~\cite{xx06} and so on, list in the next section.

Although Cayley graphs, one class of vertex-transitive graphs, are
of hight symmetry, it is not always easy to compute the distance
from a fixed vertex to all other vertices for some Cayley graphs.
For example, The $n$-dimensional cube-connected cycle $CCC_n$,
constructed from $Q_n$ by replacing each of its vertices with a
cycle $C_n=(0,1,\ldots,n-1)$ of length $n$, is a Cayley graph proved
by Carlsson {\it et al.}~\cite{ccsw85}. Until now, one has not yet
determined exactly its sum of all pairs of vertices, and so only can
give its forwarding indices asymptotically (see, Shahrokhi and
Sz\'ekely~\cite{ss01}, and Yan {\it et al}.~\cite{yxy08}).

Unfortunately, for the edge-forwarding index, there is no an analogy
of (\ref{e4}). But the lower bound of $\pi(G)$ given in (\ref{e2})
is useful. One may design a routing $R$ such that $\pi (G,R)$
attains this lower bound. For example, the edge-forwarding indices
of the folded cube~\cite{hxx05} and the augmented cube~\cite{xx06}
are determined by this method.

Making use of results on the Cartesian product is one of methods
determining forwarding indices. Using Theorem 2.15, Xu {\it et
al}~\cite{xxh06} determined the vertex-forwarding indices and the
edge-forwarding indices for the generalized hypercube
$Q(d_1,d_2,\ldots,d_n)$, the undirected toroidal mesh
$C(d_1,d_2,\ldots, d_n)$, the directed toroidal mesh
$\overrightarrow C(d_1,d_2,\ldots$, $d_n)$, all of which can be
regarded as the Cartesian products.

\vskip10pt \noindent{\bf 5.3.\quad Relations to Other Topics}

To the knowledge of the author, until now one has not yet find an
approximation algorithm with good performance ratio for finding
routings of general graphs. However, de le Vega and
Manoussakis~\cite{fm94} showed that the problem of determining the
value of the forwarding index (respectively, the forwarding index of
minimal routings) is an instance of the multicommodity flow problem
(respectively, flow with multipliers). Since many very good
heuristics or approximation algorithms are known for these flow
problems~\cite{a78, h69, l76}, it follows from these results that
all of these algorithms can be used for calculating the forwarding
index.

Teranishi~\cite{t02}, Laplacian spectra and invariants of graphs.

\vskip10pt

\section{Forwarding Indices of Some Graphs}

The forwarding indices of some very particular graphs have been
determined. We list all main results that we are interested in and
have known, all of which not noted can be found in Heydemann et
at~\cite{hms89} or determined easily.

\vskip6pt{\bf 1.}\quad For a complete graph $K_n$,\ $\xi(K_n)=0$
and $\pi(K_n)=2$.

\vskip6pt{\bf 2.}\quad For a star $K_{1,\,n-1}$,\
$\xi(K_{1,\,n-1})=(n-1)(n-2)$ and $\pi(K_{1,\,n-1})=2(n-1)$.

\vskip6pt{\bf 3.}\quad For a path $P_n$,
$\xi(P_n)=2\left\lfloor\frac
12n\right\rfloor\left(\left\lceil\frac 12n\right\rceil-1\right)$
and $\pi(P_n)=2\left\lfloor\frac 12n\right\rfloor\left\lceil\frac
12n\right\rceil$.

\vskip6pt{\bf 4.}\quad For the complete bipartite $K_{m,\,n}$
$(m\ge n)$,\ $\xi_m(K_{m,\,n})=\xi(K_{m,\,n})= \lceil\frac
{m(m-1)}n\rceil$ and $\pi_m(K_{m,1})=2m$ and if $2\le n\le m$,
$$
\left\lceil\frac {2m(m-1)+2n(n-1)}{mn}\right\rceil+2\le
\pi_m(K_{m,\,n})\le \left\lceil\frac {m-1}{n}\right\rceil.
$$
In particular,
$$
\pi_m(K_{n,\,n})=\pi_m(K_{n,\,n})=\left\{\begin{array}{ll} 4, &\
{\rm for}\ n=2;\\
5, &\ {\rm for}\ n=3,\ 4;\\ 6, &\ {\rm for}\ n\ge
5;\end{array}\right.
$$

\vskip6pt {\bf 5.}\quad For a directed cycle $C_d$ $(d\ge 3)$,\
$\xi(C_d)=\frac 12(d-1)(d-2)$. For an undirected cycle $C_d$
$(d\ge 3)$,\ $\xi_m=\xi(C_d)=\left\lfloor\frac
14(d-1)^2\right\rfloor$ and $\pi_m=\pi(C_d)=\left\lfloor\frac
14d^2\right\rfloor$.

\vskip6pt{\bf 6.}\quad The $n$-dimensional undirected toroidal mesh
$C(d_1,d_2,\cdots, d_n)$ is defined as cartesian product
$C_{d_1}\times C_{d_2}\times\cdots \times C_{d_n}$ of $n$ undirected
cycles $C_{d_1},C_{d_2},\cdots$, $C_{d_n}$ of order $d_1,d_2,\cdots,
d_n$, $d_i\ge 3$ for $i=1,2,\cdots,n$. The $C(d,d,\cdots,d)$,
denoted by $C_n(d)$, is called a $d$-ary $n$-cube a generalized
$n$-cube. Xu {\it et al}~\cite{xxh06} determined that
\begin{eqnarray*}
\xi(C(d_1,d_2,\cdots, d_n)) & = & \sum_{i=1}^{n}d_1d_2\cdots
d_{i-1}\left\lfloor \frac {d_i^2}4\right\rfloor d_{i+1}\cdots
d_n-d_1d_2\cdots d_n +1;\\
\pi(C(d_1,d_2,\cdots, d_n)) & = & \max\limits_{1\le i\le n}
\left\{d_1d_2\cdots d_{i-1}\left\lfloor\frac {d_i^2}4\right\rfloor
d_{i+1}\cdots d_n\right\}.
\end{eqnarray*}
In particular,
$$
\xi(C_n(d))=nd^{n-1}\left\lfloor \frac 14
d^2\right\rfloor-(d^n-1), \quad {\rm and}\quad
\pi(C_n(d))=d^{n-1}\left\lfloor \frac 14 d^2\right\rfloor.
$$
The last result was obtained by Heydemann et al~\cite{hms89}.

\vskip6pt{\bf 7.}\quad The $n$-dimensional directed toroidal mesh
$\overrightarrow C(d_1,d_2,\cdots, d_n)$ is defined as the cartesian
product $\overrightarrow C_{d_1}\times \overrightarrow
C_{d_2}\times\cdots \times \overrightarrow C_{d_n}$ of $n$ directed
cycles $\overrightarrow C_{d_1},\overrightarrow C_{d_2},\cdots,
\overrightarrow C_{d_n}$ of order $d_1,d_2,\cdots, d_n$, $d_i\ge 3$
for each $i=1,2,\cdots,n$. Set $\overrightarrow
C_n(d)=\overrightarrow C(d,d,\cdots,d)$. Xu {\it et al}~\cite{xxh06}
determined that
\begin{eqnarray*}
\xi(\overrightarrow{C}(d_1,d_2,\cdots, d_n)) & = &
\frac 12\left(\sum_{i=1}^n(d_i-3)\right)d_1d_2\cdots d_n+(n-1)d_1d_2\cdots d_n +1;\\
\pi(\overrightarrow{C}(d_1,d_2,\cdots, d_n)) & = & \frac 12
\max_{1\le i\le n} \{d_1\cdots d_{i-1}d_i(d_i-1)d_{i+1}\cdots
d_n\}.
\end{eqnarray*}
In particular,
$$
\xi(\overrightarrow{C}_n(d))=\frac n2\,d^n(d-1)-d^n+1,\quad {\rm
and}\quad \pi(\overrightarrow{C}_n(d))=\frac 12\,d^n(d-1).
$$

\vskip6pt {\bf 8.}\ The $n$-dimensional generalized hypercube,
denoted by $Q(d_1,d_2,\cdots,d_n)$, where $d_i\ge 2$ is an integer
for each $i=1,2,\cdots,n$, is defined as the cartesian products
$K_{d_1}\times K_{d_2}\times \cdots \times K_{d_n}$. If
$d_1=d_2=\cdots =d_n=d\ge 2$, then $Q(d,d,\cdots,d)$ is called the
$d$-ary $n$-dimensional cube, denoted by $Q_n(d)$. It is clear that
$Q_n(2)$ is $Q_n$. Xu {\it et al}~\cite{xxh06} determined that
\begin{eqnarray*}
\xi(Q(d_1,d_2,\cdots,d_n)) &=& -\sum\limits_{i=1}^n d_1d_2\cdots
d_{i-1}d_{i+1}\cdots d_n+(n-1)d_1d_2\cdots d_n+1,\\
\pi(Q(d_1,d_2,\cdots,d_n)) &=& \max\limits_{1\le i\le n}
\{d_1d_2\cdots d_{i-1}2d_{i+1}\cdots d_n\}.
\end{eqnarray*}
In particular,
$$
\xi(Q_n(d))= ((d-1)n-d)d^{n-1}+1,\quad {\rm and}\quad \pi(Q_n(d))=
2d^{n-1}.
$$
For the $n$-dimensional hypercube $Q_n$,
$$ \xi(Q_n)=(n-2)2^{n-1}+1\quad {\rm and}\quad
\pi(Q_n)=2^n.
$$
The last result was obtained by Heydemann et al~\cite{hms89}.

\vskip6pt{\bf 9.}\quad For the crossed cube $CQ_n$ $(n\ge 2)$,\
$\pi (CQ_n)=\pi_m(CQ_n)=2^n$ decided by Chang {\it et
al}~\cite{csh00}. However, $\xi (CQ_n)$ has not been determined so
far.

\vskip6pt {\bf 10.}\quad For the folded cube $FQ_n$, decided by Hou
{\it et al}~\cite{hxx05},
$$\xi(FQ_n)=\xi_m(FQ_n)=(n-1)2^{n-1}+1-\frac{n+1}{2}{n\choose
\lceil\frac n2\rceil},$$
$$\pi(FQ_n)=\pi_m(FQ_n)=2^n-{n\choose
\lceil\frac n2\rceil}.$$

\vskip6pt {\bf 11.}\quad For the augmented cube $AQ_n$ proposed by
Choudum and Sunitha~\cite{cs02}, Xu and Xu~\cite{xx06} showed that
 $$
 \xi(AQ_n)=\frac{2^n}{9}+\frac{(-1)^{n+1}}{9}+\frac{n2^n}{3}-2^n+1,
 $$
and
 $$
 \pi(AQ_n)=2^{n-1}.
 $$

\vskip6pt {\bf 12.}\quad For the cube-connected cycle $CCC(n)$ and
the $k$-dimensional wrapped butterfly $WBF_k(n)$, Yan, Xu, and
Yang~\cite{yxy08}, Shahrokhi and Sz\'ekely~\cite{ss01}, determined
 $$
\xi (CCC_n) =\frac{7}{4} n^2 2^n (1 - o(1)),
  $$
$$
\pi(CCC(n))=\pi_m (CCC(n))=\frac 54n^22^n(1-o(1)),
$$
$$
\pi(WBF(n))=\pi_m (WBF_2(n))=\frac 54n^22^{n-1}(1+o(1)).
$$
Hou, Xu and Xu~\cite{hxx09} determined
 $$
 \xi(WB_k(n))=\frac{3n(n-1)}2k^n-\frac{n(k^n-1)}{k-1}+1.
 $$

\vskip6pt {\bf 13.}\quad For the star graph $S_n$,
Gauyacq~\cite{g97} obtained that
$$
2(n-1)!(n-1)+\lceil 2\alpha\rceil\le \pi(S_n)\le
2(n-1)!(n-1)+2\lceil \alpha\rceil,
$$
where $\alpha=(n-2)!\sum\limits_{i=2}^{n-1}\frac{n-i}i$.

\vskip6pt {\bf 14.}\quad For the complete-transposition graph
$CT_n$, Gauyacq~\cite{g97} obtained that
$$
2(n-2)!(2n-3)-\lfloor 2\beta\rfloor\le \pi(CT_n)\le
2(n-2)!(2n-3)-2\lfloor \beta\rfloor,
$$
where $\beta=2(n-2)!\sum\limits_{i=3}^{n}\frac{1}i$.

\vskip6pt {\bf 15.}\quad For the undirected de Bruijn graph
$UB(d,\,n)$ and Kautz graph $UK(d,\,n)$, the upper bounds of their
vertex-forwarding indices and edge-forwarding indices have been,
respectively, given as follows.

$$\xi(UB(d,\,n))\le (n-1)d^n,\quad \xi (UK(d,\,n))\le (n-1)d^n,$$
$$\pi(UB(d,\,n))\le 2nd^{n-1},\quad \pi (UK(d,\,n))\le
2(n-1)d^{n-2}(d+1).$$

\end{document}